\documentclass{amsart} 

\newcommand{\R}{\mathbb{R}}            
\newcommand{\Z}{\mathbb{Z}}

\newcommand{\Col}{{\rm Col}}

\newcommand{\subsetm}{\stackrel{{\rm m}}{\subset}}
\newcommand{\zero}{\text{\boldmath $0$}} 

\newtheorem{theorem}{Theorem}[section]
\newtheorem{lemma}[theorem]{Lemma}
\newtheorem{corollary}[theorem]{Corollary}
\newtheorem{proposition}[theorem]{Proposition}

\newtheorem{remark}[theorem]{Remark}

\usepackage{graphics} 

\begin{document}


\title[Ribbon concordance of surface-knots]
{Ribbon concordance of surface-knots \\
via quandle cocycle invariants} 


\author{J. Scott Carter} 
\address{Department of Mathematics, 
University of South Alabama, Mobile, AL 36688, U.S.A.} 
\email{carter@mathstat.usouthal.edu} 

\author{Masahico Saito} 
\address{Department of Mathematics, 
University of South Florida, Tampa, FL 33620, U.S.A.} 
\email{saito@math.usf.edu} 

\author{Shin SATOH} 
\address{Graduate School of Science and Technology, 
Chiba University, Yayoi-cho 1-33, 
Inage-ku, Chiba, 263-8522, Japan}
\email{satoh@math.s.chiba-u.ac.jp}


\renewcommand{\thefootnote}{\fnsymbol{footnote}}





\maketitle


\begin{abstract} 
We give necessary conditions of 
a surface-knot to be ribbon concordant to another, 
by introducing a new variant of 
the cocycle invariant of surface-knots 
in addition to using the invariant already known. 
We demonstrate that twist-spins of some torus knots 
are not ribbon concordant to their orientation reversed images. 
\end{abstract}


\section{Introduction}\label{sec1}

\begin{figure}[htb]
\begin{center}
\includegraphics{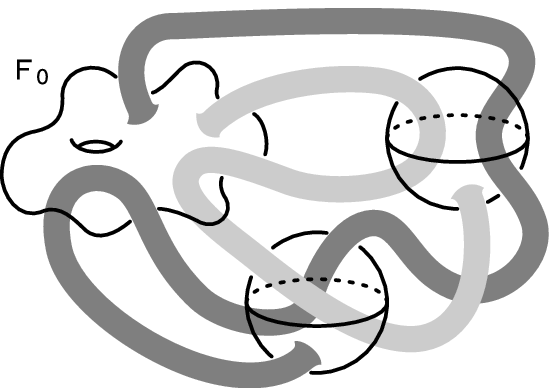}
\caption{}
\label{fig01}
\end{center}
\end{figure}

Throughout this paper, 
a {\it surface-knot} means 
a connected, oriented closed surface 
smoothly embedded 
in $4$-space ${\R}^4$ up to ambient isotopies. 
Let $F_0$ and $F_1$ be surface-knots of the same genus. 
We say that {\it $F_1$ is ribbon concordant to $F_0$} 
if there is a concordance $C$ 
in ${\R}^4\times [0,1]$ between 
$F_1\subset{\R}^4\times\{1\}$ and 
$F_0\subset{\R}^4\times\{0\}$ 
such that the restriction to $C$ 
of the projection ${\R}^4\times[0,1]\rightarrow[0,1]$ 
is a Morse function with critical points 
of index $0$ and $1$ only. 
We write $F_1\geq F_0$. 
Note that if $F_1\geq F_0$, 
then there is a set of $n$ $1$-handles 
on a split union of $F_0$ and $n$ trivial sphere-knots, 
for some $n\geq 0$, 
such that $F_1$ is obtained by surgeries 
along these handles (Figure~\ref{fig01}).

The notion of ribbon concordance was 
originally introduced by Gordon \cite{Gor} 
for classical knots in ${\R}^3$, 
and there are several studies 
found in \cite{Gil, Miy1, Miy2, Sil}, 
for example. 
Note that $F$ is a ribbon surface-knot 
if and only if 
$F$ is a ribbon concordant to the trivial sphere-knot. 

Given surface-knots $F_0$ and $F_1$, 
it is natural to ask whether 
$F_1$ is ribbon concordant to $F_0$. 
Cochran \cite{Coc} gave 
a necessary condition 
for a sphere-knot $F$ to be ribbon 
in terms the knot group $\pi_1({\R}^4\setminus F)$. 
The aim of this paper 
is to give new necessary conditions for 
a pair of surface-knots to be ribbon concordant 
by using quandle cocycle invariants. 

A quandle \cite{Joy, Mat} is an algebraic object 
whose model is a group 
with conjugation, 
and its cohomology theory 
was developed in \cite{CJKLS} 
as a generaliztion of the theory given in \cite{FRS}. 
It is known that 
each quandle $3$-cocycle $\theta$ 
defines an invariant of a surface-knot $F$, 
called the {\it quandle cocycle invariant}, $\Phi_{\theta}(F)$. 
The invariant $\Phi_{\theta}(F)$ is 
regarded as 
a multi-set of elements in the coefficient group $A$ 
of the cohomology 
where repetitions of the same element are allowed. 
For two multi-sets $A'$ and $A''$ of $A$, 
we use the notation $A'\subsetm A''$ 
if for any $a\in A'$ it holds that $a\in A''$. 
In other words, 
$A'\subsetm A''$ if and only if 
$\tilde{A'}\subset \tilde{A''}$ 
where $\tilde{A'}$ and $\tilde{A''}$ are 
the subsets of $A$ obtained from $A'$ and $A''$ 
by eliminating the multiplicity of elements, 
respectively. 
The following is a necessary condition 
for ribbon concordance.

\begin{theorem}\label{thm11} 
If $F_1\geq F_0$, 
then $\Phi_{\theta}(F_1)\subsetm\Phi_{\theta}(F_0)$. 
\end{theorem} 

By Theorem~\ref{thm11}, 
we give many examples of pairs of surface-knots 
such that one is {\it not} 
ribbon concordant to another (Corollary~\ref{cor21}). 
For example, we easily see that 
the $2$-twist-spun trefoil 
and its mirror image are not 
ribbon concordant to each other. 
However, 
Theorem~\ref{thm11} is not effective 
in the family of ribbon surface-knots; 
in fact, $\Phi_{\theta}(F)=\zero$ 
for any ribbon surface-knot $F$. 
Here, we use the notation $\zero$ 
to stand for a multi-set 
consisting of zero elements of $A$ only. 
In this paper, 
we define a new variation of cocycle invariants 
of surface-knots by using a quandle $2$-cocycle $\phi$ 
(the definition is given in Section~\ref{sec3}). 
The invariant of a surface-knot $F$ is denoted by 
$\Omega_{\phi}(F)=
\{A_{\lambda}|\lambda \in H_1(F;{\Z})\}$ 
which is a family of multi-sets $A_{\lambda}$ 
of the coefficient group $A$. 
Note that a $2$-cocycle $\phi$ is originally used 
to define the invariant, $\Phi_{\phi}(K)$, 
of a classical knot $K$ (cf. \cite{CJKLS}). 
The invariant $\Omega_{\phi}$ gives 
another necessary condition 
for ribbon concordance. 

\begin{theorem}\label{thm12} 
If $F_1\geq F_0$, 
then for any 
$A'\in\Omega_{\phi}(F_1)$, 
there is $A''\in\Omega_{\phi}(F_0)$ 
such that $A'\subsetm A''$. 
\end{theorem} 

As an application of our new invariant $\Omega_{\phi}$ 
of a surface-knot, 
we obtain a result on the cocycle invariant of 
a classical knot as follows 
(refer to \cite{Joy} for the definition of an {\it involutory} quandle, 
or see Section~\ref{sec4}). 

\begin{theorem}\label{thm13} 
If $\phi$ is a $2$-cocycle of an involutory quandle, 
then $\Phi_{\phi}(K)=\zero$ 
for any $2$-bridge knot $K$. 
\end{theorem} 

This paper is organized as follows. 
In Section~\ref{sec2}, 
we review the definition of the original cocycle invariant 
$\Phi_{\theta}(F)$. 
The proof of Theorem~\ref{thm11} 
and its application (Corollary~\ref{cor21}) 
are also contained in this section. 
In Section~\ref{sec3}, 
we introduce a new invariant $\Omega_{\phi}(F)$ 
by using a $2$-cocycle $\phi$, 
and then prove Theorem~\ref{thm12}. 
An application of the theorem is 
given in Section~\ref{sec4} (Corollary~\ref{cor43}), 
where we only sketch the outline of the proof 
and its completion is left to Appendix. 
Boyle \cite{Boy} studied a surface-knot 
obtained from a twist-spun knot 
by surgery along a $1$-handle. 
By using his result, 
we prove Theorem~\ref{thm13} also in Section~\ref{sec4}. 

\begin{remark} 
{\rm Kawauchi points out that 
the linking signature of a certain family of surface-knots 
is invariant under ribbon concordance. 
This result has not appeared in any paper, 
but can be obtained as a corollary of \cite{Kaw}. 
}
\end{remark} 


\section{Invariants by using $3$-cocycles}\label{sec2}

We first review the definition of 
the quandle $3$-cocycle invariants of surface-knots. 
Refer to \cite{CJKLS} for more details. 
A {\it quandle} 
is a set $X$ with a binary operation 
$(a,b)\mapsto a*b$ satisfying 
the following three axioms: 
\begin{itemize} 
\item $a*a=a$ for any $a\in X$, 
\item the map $*a:X\rightarrow X$ defined by 
$x\mapsto x*a$ is bijective for any $a\in X$, 
and 
\item $(a*b)*c=(a*c)*(b*c)$ 
for any $a,b,c\in X$. 
\end{itemize} 
For an abelian group $A$, 
we say that a map 
$\theta:X^3\rightarrow A$ 
is a {\it $3$-cocycle} if it satisfies the conditions that 
\begin{itemize} 
\item 
$\theta(x_1,x_2,x_3)=0$ 
if $x_1=x_2$ or $x_2=x_3$, and 
\item 
$\theta(x_1,x_3,x_4)-\theta(x_1,x_2,x_4)+\theta(x_1,x_2,x_3)$

\noindent 
$=\theta(x_1*x_2,x_3,x_4) 
-\theta(x_1*x_3,x_2*x_3,x_4)
+\theta(x_1*x_4,x_2*x_4,x_3*x_4)$ 

\noindent
for any $x_1,\dots,x_4\in X$. 
\end{itemize} 
We denote by $Z^3(X;A)$ 
the set of such $3$-cocycles.

To describe a surface-knot, 
we use a fixed projection of $\pi: {\R}^4\rightarrow{\R}^3$ 
as well as a description of a classical knot 
into the plane. 
Every surface-knot $F$ can be perturbed 
slightly in ${\R}^4$ 
so that the projection 
image $\pi(F)$ has double point curves, 
isolated triple points, 
and isolated branch points 
as the closures of the multiple point set. 
Crossing information is indicated in $\pi(F)$ 
as follows: 
Along every double point curve, 
two sheets intersect locally, 
one of which is under the other 
relative to the projection direction of $\pi$. 
Then the under-sheet is broken by the over-sheet. 
A {\it diagram} of $F$ 
is the image $\pi(F)$ with such crossing information. 
Hence, a diagram is regarded as a union of 
disjoint compact, connected surfaces. 
For a diagram $D$, 
we denote by $\Sigma(D)$ 
the set of such connected suraces of $D$. 
Note that three sheets near a triple point 
are labeled top, middle, and bottom 
according to crossing information, 
and the middle and bottom sheets are 
divided into two and four pieces, respectively.

\begin{figure}[htb]
\begin{center}
\includegraphics{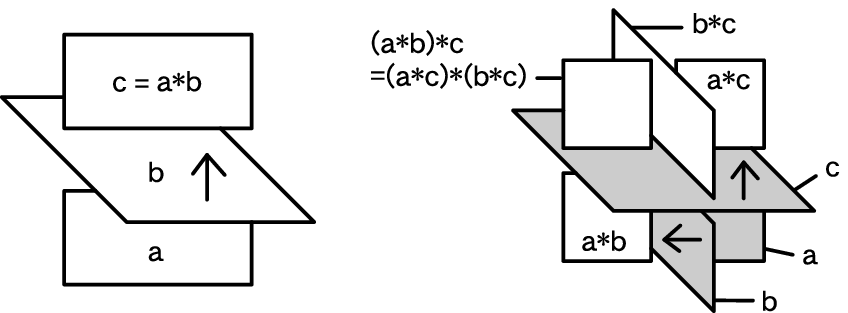}
\caption{}
\label{fig02}
\end{center}
\end{figure}

For a quandle $X$, 
a map $C:\Sigma(D)\rightarrow X$ is called 
an {\it $X$-coloring} of $D$ 
if it satisfies the following condition 
near every double point $d$: 
if $a=C(\alpha_1)$ and $c=C(\alpha_2)$ 
are the colors of under-sheets $\alpha_1$ and $\alpha_2$ 
separated by the over-sheet $\beta$ 
colored by $b=C(\beta)$, 
where the orientation normal of $\beta$ 
points from $\alpha_1$ to $\alpha_2$, 
then $a*b=c$ holds. 
See the left of Figure~\ref{fig02}. 
We denote the set of such $X$-colorings of $D$ 
by ${\Col}_X(D)$. 
Also, the pair $(a,b)$ is called the {\it color} of 
a double point $d$, 
and denoted by $C(d)\in X^2$. 

Each triple point $t$ of $D$ is assigned the sign 
$\varepsilon(t)=\pm 1$ induced from the orientation 
in such a way that $\varepsilon(t)=+1$ if and only if 
the ordered triple of the orientation normals of 
the top, middle, and bottom sheets, respectively, 
agrees with the orientation of ${\R}^3$. 
Given an $X$-coloring $C\in{\Col}_X(D)$, 
the colors of the sheets near $t$ are 
determined by three colors 
$a=C(\alpha)$, $b=C(\beta)$, and $c=C(\gamma)$, 
where $\gamma$ is the top sheet, 
$\beta$ is the middle sheet from which 
the orientation normal of $\gamma$ points, 
and $\alpha$ is the bottom sheet 
from which the orientation normals of 
$\beta$ and $\gamma$ point both. 
See the right of Figure~\ref{fig02}, 
where the sheets $\alpha,\beta$, and $\gamma$ are 
shaded. 
The ordered triple $(a,b,c)$ is 
called the {\it color} of $t$ and 
denoted by $C(t)\in X^3$.

Let $D$ be a diagram of $F$ 
colored by $C\in{\Col}_X(D)$. 
Given a $3$-cocycle $\theta \in Z^3(X;A)$, 
we define the ({\it Boltzmann}) {\it weight} 
of each triple point $t$ by 
$W_{\theta}(t;C)=\varepsilon(t)\cdot\theta(a,b,c)\in A$, 
where $C(t)=(a,b,c)$. 
We denote by $W_{\theta}(C)\in A$ 
the sum $\sum_{t}W_{\theta}(t;C)$ 
for all triple points of $D$. 
Then the {\it cocycle invariant} of $F$ 
by using $\theta$ is the multi-set 
$$\Phi_{\theta}(F)=
\bigl\{W_{\theta}(C)\in A|C\in{\Col}_X(D)\bigr\},$$
where repetitions of the same element are allowed. 
It is proved in \cite{CJKLS} 
to be an invariant of $F$ 
which does not depend on the choice of 
a diagram $D$ of $F$.

Let $F_0$ and $F_1$ be surface-knots 
with $F_1\geq F_0$, 
that is, $F_1$ is ribbon cocordant to $F_0$. 
For a diagram $D_0$ of $F_0$, 
we may take a typical diagram $D_1$ of $F_1$ as follows: 
There is a set of sufficiently thin 
$n$ $1$-handles, $h_1,\dots,h_n$, 
for some $n\geq 0$, 
connecting a split union of $D_0$ and 
$n$ embedded $2$-spheres, 
$S_1,\dots,S_n$, such that 
\begin{itemize} 
\item each $1$-handle $h_j$ 
connects $D_0$ and $S_j$, 
and intersects 
$D_0\cup\bigl(\bigcup_{i=1}^n S_i\bigr)$ 
with disjoint meridian $2$-disks of $h_j$, and 
\item $D_1$ is obtained from 
$D_0\cup\bigl(\bigcup_{i=1}^n S_i\bigr)$ 
by surgeries along $\bigcup_{j=1}^n h_j$. 
\end{itemize} 
In the following, 
we use $D_1$ in the above form unless otherwise stated.

\begin{proof}[Proof of Theorem~$\ref{thm11}$.] 
For any element $a\in \Phi_{\theta}(F_1)$, 
there is an $X$-coloring $C_1\in{\Col}_X(D_1)$ 
with $a=W_{\theta}(C_1)=\sum_{t}W_{\theta}(t;C_1)$ 
on $D_1$. 
Since the intersection of $D_0$ 
and each $1$-handle $h_j$ 
consists of small $2$-disks, 
the $X$-coloring $C_1$ restricted to the punctured diagram 
$D_0\setminus \bigl(\bigcup_{j=1}^n h_j\bigr)$ 
determines the $X$-coloring of $D_0$ uniquely, 
$C_0\in{\Col}_X(D_0)$. 
Since the set of triple points of $D_1$ is 
coincident with that of $D_0$, 
and since $W_{\theta}(t;C_0)=W_{\theta}(t;C_1)$ 
for any triple point $t$, 
we have 
$a=\sum_{t}W_{\theta}(t;C_0)=
W_{\theta}(C_0) \in \Phi_{\theta}(D_0)$. 
\end{proof}

We present specific examples 
as an application of 
Theorem~\ref{thm11} in the rest of this section. 
The set $\{0,1,\dots,p-1\}$ becomes a quandle 
under the operation 
$a*b=2b-a$ (mod $p$), 
which is called the {\it dihedral quandle} of order $p$, 
and denoted by $R_p$. 
For an odd prime $p$, 
Mochizuki \cite{Moc} found 
a $3$-cocycle $\theta_p\in Z^3(R_p,{\Z}_p)$ given by 
$$\theta_p(x_1,x_2,x_3)=
(x_1-x_2)\frac{(2x_3-x_2)^p+x_2^p-2x_3^p}{p},$$
where coefficients in the numerator are divisible by $p$. 
The reader can check that $\theta_p$ 
satisfies the $3$-cocycle conditions by hands (cf. \cite{AS}).

In $1965$ Zeeman \cite{Zee} introduced 
an important family of sphere-knots. 
We take a tangle (knotted arc) $T_K$ in the $3$-ball $B^3$, 
whose closure is a classical knot $K$. 
For an integer $r\geq 0$, 
let $\{f_t\}_{t\in[0.1]}$ be the ambient isotopy of $B^3$ 
which rotates the tangle $T_K$ 
a total of $r$ times about an axis 
while keeping the boundary of $T_K$ fixed. 
Futhermore, $f_0(T_K)=f_1(T_K)$. 
We construct an annulus $A$ 
properly embedded in $B^3\times S^1$ from 
$$\bigcup_{t\in[0,1]} f_t(T_K) \times \{t\} 
\subset B^3\times [0,1]$$
by identifying the quotient $[0,1]/(0=1)$ with $S^1$. 
The {\it $r$-twist-spin of $K$} is a sphere-knot obtained 
by embedding $(B^3\times S^1, A)$ in ${\R}^4$ standardly 
and capping $A$ with two $2$-disks 
along the boundary of $A$. 
We denote the sphere-knot by $\tau^rK$. 

Let $T(2,q)$ denote the $(2,q)$-torus knot in ${\R}^3$. 
For a surface-knot $F$, 
let $-F$ denote the surface-knot $F$ 
with the reversed orientation. 
Then we have the following.

\begin{corollary}\label{cor21} 
{\rm (i)} If $q$ and $q'$ are distinct odd primes, 
then we have 
$\tau^2T(2,q)\not\geq \tau^2T(2,q')$ and 
$\tau^2T(2,q')\not\geq \tau^2T(2,q)$. 

{\rm (ii)} If $q$ is an odd prime with 
$q\equiv 3$ {\rm (mod $4$)}, 
then we have 
$\tau^2T(2,q)\not\geq -\tau^2T(2,q)$ 
and $-\tau^2T(2,q)\not\geq \tau^2T(2,q)$. 
\end{corollary} 


\begin{proof} 
(i) It is proved in \cite{AS} that 
$\Phi_{\theta_p}\bigl(\tau^2T(2,q)\bigr)=\zero$ 
for $p\ne q$, and 
$$\renewcommand{\arraystretch}{1.3} 
\Phi_{\theta_q}\bigl(\tau^2T(2,q)\bigr)=\left\{
\begin{array}{rrr}
0, & \dots, & 0, \\
-2\cdot 1^2, & \dots, & -2\cdot 1^2, \\
-2\cdot 2^2, & \dots, & -2\cdot 2^2, \\
\dots & \dots & \dots, \\
-2\cdot (q-1)^2, & \dots, & -2\cdot (q-1)^2 
\end{array}\right\}$$
for $p=q$, 
where the number of each term 
of the form $-2k^2$ $(k=0,1,\dots,q-1)$ is $q$. 
In particular, 
since 
$\Phi_{\theta_q}\bigl(\tau^2T(2,q)\bigr)\ne\zero$ 
and $\Phi_{\theta_q}\bigl(\tau^2T(2,q')\bigr)=\zero$, 
we have 
$\tau^2T(2,q)\not\geq \tau^2T(2,q')$ by Theorem~\ref{thm11}. 
It is also similarly proved that 
$\tau^2T(2,q')\not\geq \tau^2T(2,q)$. 

(ii) It is known that 
$\Phi_{\theta}(-F)=-\Phi_{\theta}(F)$ 
for any surface-knot $F$ and $3$-cocycle $\theta$ 
(cf.~\cite{CJKLS}). 
On the other hand, 
we obtain the set 
$S=\bigl\{-2k^2|k=0,1,\dots,\frac{p-1}{2}\bigr\}$ 
from $\Phi_{\theta_q}\bigl(\tau^2T(2,q)\bigr)$ 
by eliminating the multiplicity of elements. 
It is not difficult to see that if $q\equiv 3$ (mod $4$), 
then $S\not\subset -S$ and $S\not\supset -S$, 
and hence, we have the conclusion 
by Theorem~\ref{thm11}. 
\end{proof}


\section{Invariants by using $2$-cocycles}\label{sec3}

Let $X$ be a quandle and $A$ an abelian group. 
We say that a map 
$\phi:X^2\rightarrow A$ 
is a {\it $2$-cocycle} if it satisfies 
\begin{itemize} 
\item 
$\phi(x_1,x_2)=0$ 
if $x_1=x_2$, and 
\item 
$\phi(x_1,x_3)-\phi(x_1,x_2)
=\phi(x_1*x_2,x_3) -\phi(x_1*x_3,x_2*x_3)$ 
for any $x_i\in X$. 
\end{itemize} 
We denote by $Z^2(X;A)$ 
the set of such $2$-cocycles. 

\begin{figure}[htb]
\begin{center}
\includegraphics{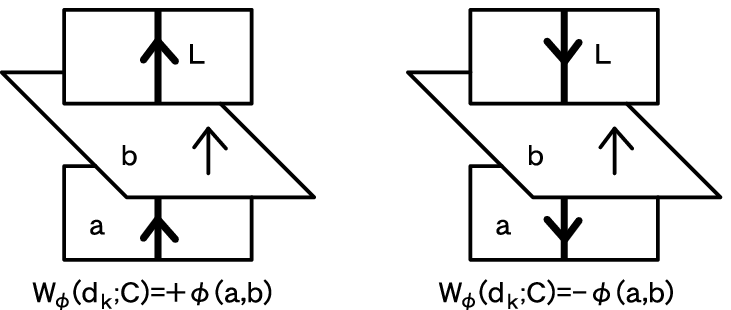}
\caption{}
\label{fig03}
\end{center}
\end{figure}

We define a cocycle invariant of a surface-knot 
by using a $2$-cocycle $\phi\in Z^2(X;A)$. 
Let $D$ be a diagram of a surface-knot $F$, 
and $C\in{\Col}_X(D)$ an $X$-coloring of $D$. 
Consider an oriented immersed circle $L$ on $D$ 
intersecting the double point curves transversely, 
and missing triple points and branch points. 
Let $d_1,\dots,d_m$ denote the points 
on the under-sheet at which 
$L$ intersects the double point curves. 
We give the sign $\varepsilon(d_k)=\pm 1$ to $d_k$ 
such that $\varepsilon(d_k)=+1$ if and only if 
the orientation of $L$ at $d_k$ agrees with 
the orientation normal of the over-sheet. 
We define the Boltzman weight at $d_k$ 
by $W_\phi(d_k;C)=\varepsilon(d_k)\cdot \phi(a,b)\in A$, 
where $C(d_k)=(a,b)$. 
Moreover, 
we put $W_{\phi}(L;C)=
\sum_{k=1}^m W_{\phi}(d_k;C)$. 
See Figure~\ref{fig03}. 
We extend these notations 
for a union of immersed circles $L$ on $D$ 
naturally. 

\begin{lemma}\label{lem31} 
If $L$ and $L'$ are homologous on $D$, 
then $W_{\phi}(L;C)=W_{\phi}(L';C)$. 
\end{lemma}

\begin{figure}[htb]
\begin{center}
\includegraphics{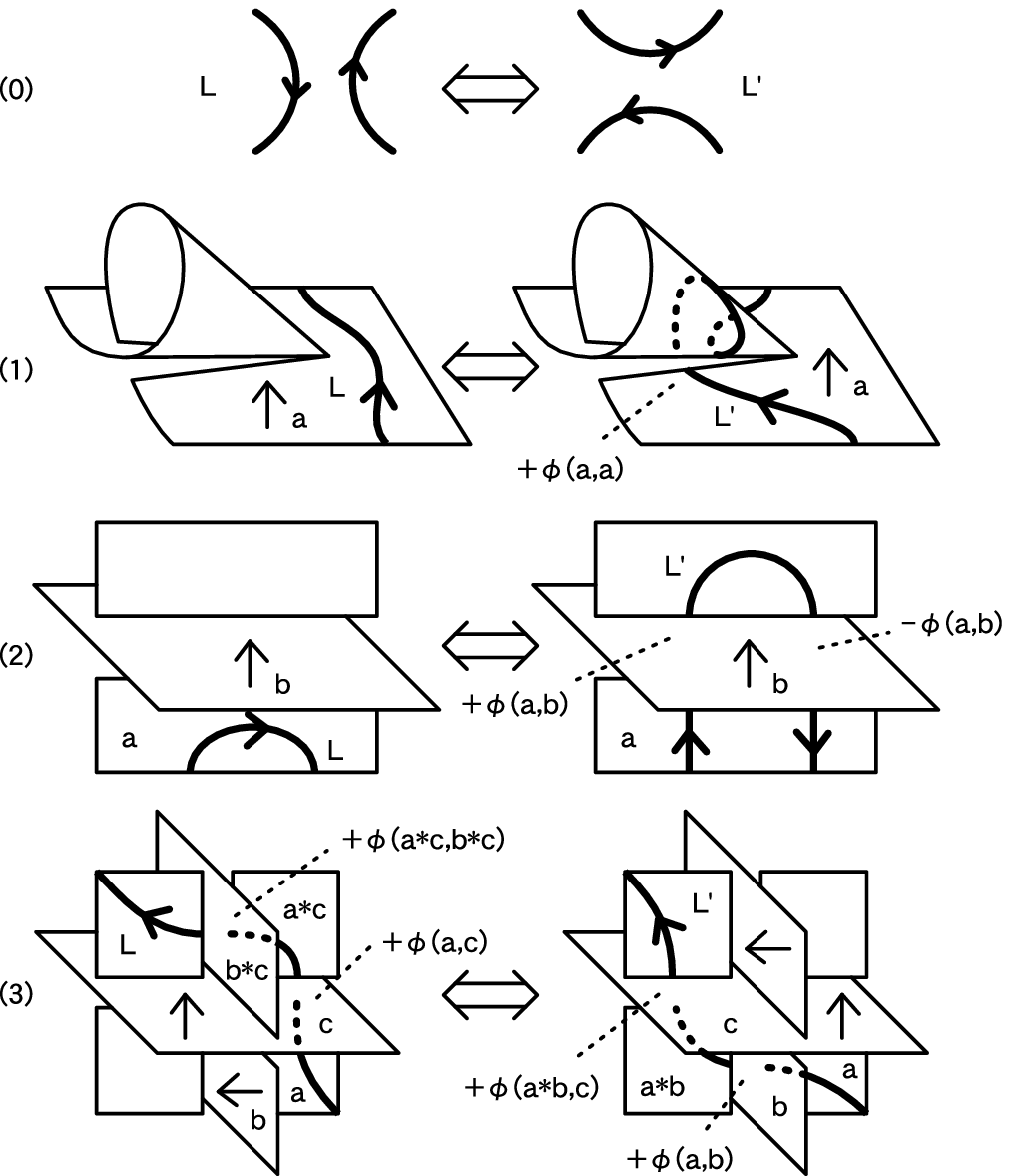}
\caption{}
\label{fig04}
\end{center}
\end{figure}

\begin{proof} 
It is sufficient to prove that 
$W_{\phi}(L;C)$ does not change 
under the moves $(0)$--$(3)$ 
(and the ones with orientation reversed, 
or with opposite crossing information) as shown in 
Figure~\ref{fig04}. 
First, 
it is clear for the move $(0)$ by the definition of 
$W_{\phi}(L;C)$. 
Since $\phi$ satisfies 
$\phi(a,a)=0$ for any $a\in X$, 
the move $(1)$ also does not change $W_{\phi}(L;C)$. 
In the move $(2)$, 
the terms $\phi(a,b)$ and $-\phi(a,b)$ 
cancels in $W_{\phi}(L';C)$. 
Finally, 
it follows from the $2$-cocycle condition of $\phi$ 
that $W_{\phi}(L;C)=W_{\phi}(L';C)$ 
under the move $(3)$. 
\end{proof} 

For each homology class $\lambda\in H_1(F;{\Z})$ 
and its representative curve $L\subset D$, 
the element $W_{\phi}(L;C)\in A$ 
is independent of the choice of $L$ by Lemma~\ref{lem31}, 
and hence, 
we denote it by $W_{\phi}(\lambda;C)$. 
Then we assign each class $\lambda\in H_1(F;{\Z})$ 
a multi-set $\Omega_{\phi}(\lambda)$ of $A$ such that 
$$\Omega_{\phi}(\lambda)=
\bigl\{W_{\phi}(\lambda;C)\bigl| C\in{\Col}_X(D)\bigr\}.$$ 
Moreover, 
we define a family of multi-sets of $A$ by 
$$\Omega_{\phi}(F)=
\bigl\{\Omega_{\phi}(\lambda)\bigl| 
\lambda\in H_1(F;{\Z})\bigr\}. $$

\begin{proposition}\label{prop32} 
The family $\Omega_{\phi}(F)$ 
does not depend on the choice of 
a diagram $D$ of $F$. 
\end{proposition} 

\begin{proof} 
It is known that any other diagram $D'$ of $F$ is 
obtained from $D$ by a finite sequence of 
Roseman moves \cite{Ros} 
up to ambient isotopies of ${\R}^3$. 
Assume that $D'$ is obtained from $D$ 
by a single Roseman move 
in a sufficiently small $3$-ball $B^3$. 
For any class $\lambda\in H_1(F;{\Z})$, 
we may take its representative curve $L$ on $D$ 
with $L\cap B^3=\emptyset$ 
so that we regard $L$ as a curve on $D'$ also. 
Moreover, 
each $X$-coloring $C\in{\Col}_X(D)$ induces 
a coloring $C'\in{\Col}_X(D')$ uniquely. 
Hence, 
any $W_{\phi}(L;C)$ on $D$ 
is coincident with 
$W_{\phi}(L;C')$ on $D'$. 
\end{proof}

The following proof is similar to 
that of Theorem~\ref{thm11} 
in Section~\ref{sec2}.

\begin{proof}[Proof of Theorem~$\ref{thm12}$.] 
Let $D_i$ be diagrams of $F_i$ $(i=0,1)$ 
as in the beginning of Section~\ref{sec3}. 
For any $A'\in\Omega_{\phi}(F_1)$, 
there is a curve $L$ on $D_1$ with 
$A'=\Omega_{\phi}(L)$. 
Since we can deform $L$ 
such that $L\cap\bigl(\bigcup_{j=1}^n h_j\bigr)=\emptyset$, 
$L$ is regarded as a curve on $D_0$. 
Put $A''=\Omega_{\phi}(L)\in\Omega_{\phi}(F_0)$. 
Then $A'\subsetm A''$ can be proved in a similar way to 
Theorem~\ref{thm11}. 
\end{proof}


\section{Torus-knots with $1$-handles}\label{sec4} 

A surface-knot is called a {\it torus-knot} 
if it is an embedded torus in ${\R}^4$. 
We distinguish it from a classical `torus knot' in ${\R}^3$ 
by inserting the hyphen -. 
In this section, 
we use a typical family of 
torus-knots studied by Boyle \cite{Boy}. 
Let $K$ be a classical knot in a $3$-ball $B^3$, 
and let $D^3\subset{\rm int}B^3$ 
be a $3$-ball such that $D^3\cap K=T_K$ 
is the knotting arc for $K$. 
For an integer $r\geq 0$, 
let $\{g_t\}_{t\in[0,1]}$ be the ambient isotopy of $B^3$ 
which rotates $T_K$ $r$ times 
keeping the trivial arc $K \setminus T_K$ fixed. 
We denote by $\sigma^rK$ 
the torus-knot obtained from 
$\bigcup_{t}g_t(K)\times\{t\}\subset B^3\times S^1$ 
by embedding it in ${\R}^4$ standardly. 
Note that $\sigma^rK$ is also obtained from 
the $r$-twist-spin of $K$ by surgery 
along a certain $1$-handle $h$. 

\begin{figure}[htb]
\begin{center}
\includegraphics{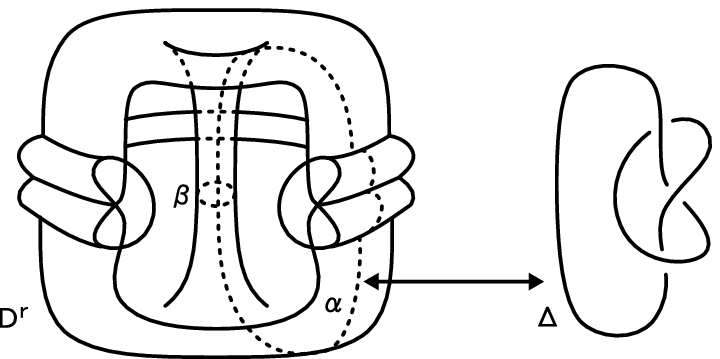}
\caption{}
\label{fig05}
\end{center}
\end{figure}

By definition, 
$\sigma^rK$ has a diagram $D^r$ 
in the form $\Delta\times S^1$, 
where $\Delta$ is a knot diagram of $K$, 
except the twisting part of $T_K$. 
See Figure~\ref{fig05}, 
where we ignore crossing information 
along double point curves and 
omit the twisting part. 
Refer to \cite{AS, Sat} for the complete 
figure of the diagram. 
We take meridional and longitudinal curves 
$\alpha$ and $\beta$ on $D^r$, respectively, 
such that $\alpha$ can be identified with $\Delta$, 
and $\beta$ has no intersection with 
double point curves.

To calculate the invariant $\Omega_{\phi}$ of $\sigma^rK$, 
we recall the definition of the cocycle invariant 
$\Phi_{\phi}(K)$ 
of a classical knot $K$ by using a $2$-cocycle $\phi$. 
Let $\Delta\subset{\R}^2$ be a diagram 
of an oriented classical knot $K$, 
and $\Sigma(\Delta)$ the set of arcs separated by 
over-arcs at crossings. 
For a quandle $X$, 
a map $C:\Sigma(\Delta)\rightarrow X$ is 
called an $X$-coloring of $\Delta$ 
if it satisfies the following condition 
near every crossing $x$: 
if $a=C(\alpha_1)$ and $c=C(\alpha_2)$ are 
the colors of under-arcs $\alpha_1$ and $\alpha_2$ 
separated by the over-arc $\beta$ colored by 
$b=C(\beta)$, 
where $\alpha_1$ is on the right side of $\beta$, 
then $a*b=c$ holds. 
We denote the set of such $X$-colorings of $\Delta$ 
by ${\Col}_X(\Delta)$. 
Also, the pair $(a,b)$ is called the color of 
the crossing $x$, and 
denoted by $C(x)\in X^2$. 
Given a $2$-cocycle $\phi\in Z^2(X;A)$, 
we define the Boltzmann wieight at $x$ by 
$W_{\phi}(x;C)=
\varepsilon(x)\cdot \phi(a,b)\in A$, 
where $C(x)=(a,b)$. 
We denote by $W_{\phi}(C)\in A$ 
the sum $\sum_{x}W_{\phi}(x;C)$ 
for all crossings of $\Delta$. 
Then the cocycle invariant of $K$ 
by using $\phi$ is the multi-set 
$$\Phi_{\phi}(K)=
\bigl\{W_{\phi}(C)\bigl| C\in{\Col}_X(\Delta)\bigr\}$$ 
where repetitions of the same element are allowed. 
It is proved in \cite{CJKLS} to be 
an invariant of $K$ which does not depend on 
the choice of a diagram $\Delta$ of $K$.

Any $X$-coloring of $D^r$ 
determines that of $\Delta$ 
by restricting it to the meridional curve $\alpha$. 
Conversely, 
not any $X$-coloring of $\Delta$ 
extend to $D^r$ totally; 
an $X$-coloring of $\Delta$ 
extends to $D^r$ 
if and only if $x(*y)^r=x$ for any 
$x,y \in X$ appeared in $\Delta$; 
this condition corresponds to the $r$-twisting of $T_K$. 
Refer to \cite{AS,Sat} for more details. 
A quandle $X$ is called {\it of type $s$} $(s\geq 0)$ 
if it satisfies that $x(*y)^s=x$ for any $x,y\in X$, 
and in particular, 
$X$ is an {\it involutory} quandle 
if it is of type $2$. 
The dihedral quandle $R_p$ is an example of 
involutory quandles; 
$(x*y)*y\equiv (2y-x)*y\equiv 2y-(2y-x)\equiv x$ 
(mod $p$). 
Then we have the following immediately.

\begin{lemma}[cf. \cite{AS,Sat}]\label{lem41} 
If $X$ is a quandle of type $s$, 
then for any $r=0, s, 2s, 3s, \dots$, 
there is a natural one-to-one correspondence 
between $\Col_X(D^r)$ and $\Col_X(\Delta)$. 
\qed 
\end{lemma} 

\begin{proposition}\label{prop42} 
Assume that $X$ is a quandle of type $s$, 
and let $\phi\in Z^2(X;A)$ a $2$-cocycle of $X$. 
For any $r=0, s, 2s, 3s, \dots$, 
the cocycle invariant $\Omega_{\phi}(\sigma^rK)$ 
is given by 
$$\renewcommand{\arraystretch}{1.3} 
\Omega_{\phi}(\sigma^rK)=
\left\{
\begin{array}{rrr}
\dots & \dots & \dots \\
-\Phi_{\phi}(K), & -\Phi_{\phi}(K), & \dots \\
\zero, & \zero, & \dots \\
\Phi_{\phi}(K), & \Phi_{\phi}(K), & \dots \\ 
2\Phi_{\phi}(K), & 2\Phi_{\phi}(K), & \dots \\
\dots & \dots & \dots 
\end{array}\right\},$$
where the number of each multi-set $k\Phi_{\phi}(K)$ 
is infinite $(k\in{\Z})$. 
In particular, 
we have $\Phi_{\phi}(K)\in\Omega_{\phi}(\sigma^rK)$. 
\end{proposition} 

\begin{proof} 
Recall that $(\alpha,\beta)$ represents 
a basis of $H_1(\sigma^rK;{\Z})\cong{\Z}\oplus{\Z}$. 
For any class $\lambda=k[\alpha]+l[\beta]$ 
$(k,l\in{\Z})$, 
we have 
$$W_{\phi}(\lambda;C)=
kW_{\phi}(\alpha;C)+lW_{\phi}(\beta;C)
=kW_{\phi}(\alpha;C)$$
by definition. 
Hence, it follows from Lemma~\ref{lem41} that 
$W_{\phi}(\lambda)=k\Phi_{\phi}(K)$. 
\end{proof}

\begin{figure}[htb]
\begin{center}
\includegraphics{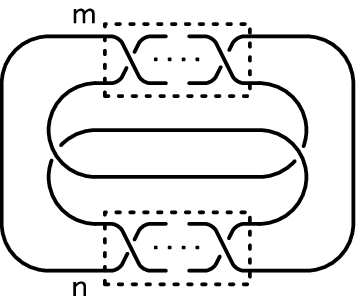}
\caption{}
\label{fig06}
\end{center}
\end{figure}

For integers $m$ and $n$, 
let $S(m,n)$ be the classical knot represented 
by the diagram as shown in Figure~\ref{fig06}. 
Note that $S(3,3)$ is coincident with 
$8_5$ in the knot table. 
Then we have the following as 
a corollary of Theorem~\ref{thm12}. 
We sketch the outline of the proof here, 
and a complete proof is given in Appendix.

\begin{corollary}\label{cor43} 
We have 
$\sigma^rT(2,l) \not\geq \sigma^sS(m,n)$ 
for any $r,s\equiv 0$ {\rm (mod~4)} and 
$l,m,n\equiv 3$ {\rm (mod 6)}. 
\end{corollary} 

\begin{proof} 
Let $Q_6$ be the subset 
of the permutation group of four letters, 
consisting of six cyclic elements of length four. 
Then $Q_6$ has a quandle structure under conjugation. 
Note that $Q_6$ is a quandle of type $4$. 
There exists a $2$-cocycle $\phi\in Z^2(Q_6;{\Z}_4)$ 
with the coefficient group ${\Z}_4$ 
such that 
the associated invariant of the $(2,k)$-torus knot 
satisfies 
$$\Phi_{\phi}\bigl(T(2,l)\bigr)=
\bigl\{\underbrace{0,0,\dots,0}_{6}, 
\underbrace{l+2,l+2,\dots,l+2}_{24}\bigr\}$$
for any $l\equiv 3$ (mod $6$), 
where the values in the invariant are taken in ${\Z}_4$. 
On the other hand, 
the invariant of $S(m,n)$ 
associated with the same $2$-cocycle $\phi$ 
satisfies 
$$\Phi_{\phi}\bigl(S(m,n)\bigr)=
\bigl\{\underbrace{0,0,\dots,0}_{30}, 
\underbrace{2,2,\dots,2}_{24}\bigr\}$$
for any $m,n\equiv 3$ (mod $6$). 
By Proposition~\ref{prop42}, 
we have $\Phi_{\phi}\bigl(T(2,l)\bigr)\in 
\Omega_{\phi}(\sigma^rT(2,l)\bigr)$ 
for $r\equiv 0$ (mod $4$). 
Since 
$$\Phi_{\phi}\bigl(T(2,l)\bigr)
\not\subsetm k\Phi_{\phi}\bigl(S(m,n)\bigr)
=\{0,0,\dots,2k,2k,\dots\}\in 
\Omega_{\phi}\bigl(\sigma^sS(m,n)\bigr)$$
for any $s\equiv 0$ (mod $4$) and $k\in {\Z}$, 
it follows from Theorem~\ref{thm12} that 
$\sigma^rT(2,l)$ is not ribbon concordant 
to $\sigma^sS(m,n)$. 
\end{proof}

To prove Theorem~\ref{thm13}, 
we prepare the following lemma. 
We say that a torus-knot 
is {\it reducible} if 
it is obtained from a sphere-knot by surgery 
along a trivial $1$-handle. 

\begin{lemma}\label{lem44} 
If a torus-knot $F$ is reducible, 
then $\Omega_{\phi}(F)=\{\zero,\zero,\dots,\zero, \dots\}$ 
for any $2$-cocycle $\phi$. 
\end{lemma} 

\begin{proof} 
For any class $\lambda\in H_1(F;{\Z})$, 
we can choose a representative curve $L$ of $\lambda$ 
along the trivial $1$-handle 
which does not meet any double point curves. 
Hence, we have $W_{\phi}(\lambda;C)=0$ by definition. 
\end{proof} 

\begin{proof}[Proof of Theorem~$\ref{thm13}$.] 
Consider the invariant $\Omega_{\phi}(\sigma^2K)$ 
of the torus-knot $\sigma^2K$. 
Since $X$ is an involutory quandle, 
that is, of type $2$, 
we have $\Phi_{\phi}(K)\in\Omega_{\phi}(\sigma^2K)$ 
by Proposition~\ref{prop42}. 
On the other hand, 
Boyle \cite{Boy} proved that if $K$ is a $2$-bridge knot, 
then $\sigma^2K$ is a reducible torus-knot. 
Hence, we have $\Phi_{\phi}(K)=\zero$ 
by Lemma~\ref{lem44}. 
\end{proof}


\section*{Achknowledgments}

The first, second, and third authors are 
partially supported by 
NSF Grant DMS $\#0301095$, 
NSF Grant DMS $\#0301089$, 
and JSPS Postdoctoral Fellowships for Research Abroad, 
respectively. 
The third author expresses his gratitude 
for the hospitality of University of South Florida.



\section*{Appendix}\label{app}

Let $Q_6$ be the subset of the permutation group 
of four letters $1,2,3,4$ 
consisting of cyclic elements of length four, 
where the subscript $6$ stands for the number of elements 
belonging to $Q_6$. 
Then $Q_6$ becomes a quandle 
under conjugation $g*h=h^{-1}gh$; 
in general, any conjugacy class of a group 
becomes a quandle under the conjugation. 
For example, if $g=(1342)$ and $h=(1234)$, 
then 
$$g*h=(1342)*(1234)=(1234)^{-1}(1342)(1234) 
=(1324).$$ 
In other words, 
$g*h$ is obtained from $g$ 
by replacing the letters in $g$ according to 
the permutation of $h$. 
Note that since $g(*h)^4=h^{-4}gh^4=g$, 
$Q_6$ is a quandle of type $4$.

\begin{figure}[htb]
\begin{center}
\includegraphics{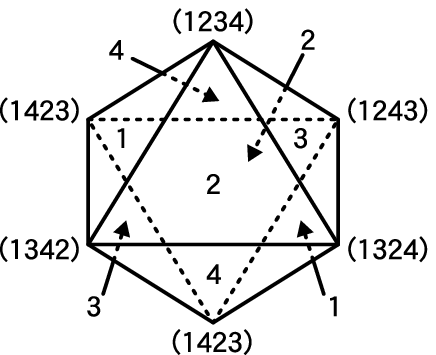}
\caption{}
\label{fig07}
\end{center}
\end{figure}

The quandle $Q_6$ can be visualized 
by using the equilateral octahedron $H$ 
(Figure~\ref{fig07}). 
First, we number the faces of $H$ 
by $1,\dots,4$ in such a way that 
each pair of parallel faces admit the same number. 
At each vertex, we put the element of $Q_6$ 
by reading the numbers on faces concentrated at the vertex 
counterclockwise. 
Under the identification of $Q_6$ and 
the set of vertices of $H$, 
the vertex $g*h$ is obtained from $g$ 
by rotating $H$ quarterly around the diagonal axis through $h$ 
in the counterclockwise direction; 
in fact, the permutation of the numbers 
on faces caused by the rotation 
is coincident with $h$ as an element of $Q_6$. 
Note that for each vertex $g\in Q_6$, 
the inverse $g^{-1}$ is located 
on the diagonal vertex. 
Quandles consisting of rotations of 
an equilateral polyhedron can be found in \cite{AG}.

Recall the the $2$-cocycle conditions are 
\begin{eqnarray} 
\phi(a,a)=0 & \mbox{for any } a\in Q_6, \mbox{ and}\\ 
\phi(a,c)-\phi(a,b)-\phi(a*b,c)+\phi(a*c,b*c)=0 & 
\mbox{for any } a,b,c\in Q_6. 
\end{eqnarray} 
By using the model of the octahedron $H$, 
we will give a way to check 
whether a given map 
$\phi:Q_6\times Q_6\rightarrow A$ 
satisfies the condition (2). 
For this purpose, 
we interpret the condition (2) visually.

\begin{figure}[htb]
\begin{center}
\includegraphics{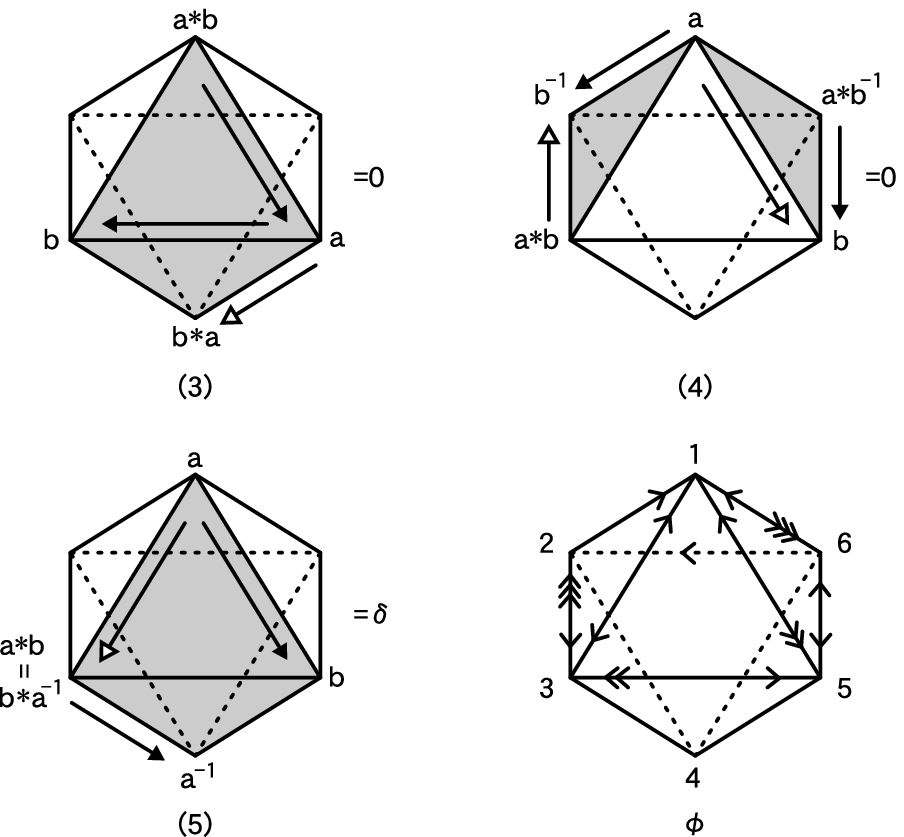}
\caption{}
\label{fig08}
\end{center}
\end{figure}

\underline{Case 1.} 
Assume that the set $\{a,b,c\}$ contains 
the same element. 
\begin{itemize} 
\item[\underline{1-i.}] If $a=b$ or $b=c$, 
then the condition(2) always holds under (1). 
\item[\underline{1-ii.}] 
Assume that $a=c\ne b$. 
If $b=a^{-1}$, 
then (2) always holds similarly. 
If $b\ne a^{-1}$, 
then (2) is equivalent to 
\begin{eqnarray} 
\phi(a,b)+\phi(a*b,a)-\phi(a,b*a)=0, 
\end{eqnarray} 
for any pair $(a,b)$ which spans an edge of 
the octahedron $H$. 
We illustrate this condition (3) 
as in Figure~\ref{fig08}, 
where the black/white arrow $\overrightarrow{xy}$ 
corresponds to the value 
$\phi(x,y)$ or $-\phi(x,y)$, respectively. 
\end{itemize}

\underline{Case 2.} 
Assume that $\{a,b,c\}$ contains 
no pair of the same element but 
a pair of inverse elements. 
\begin{itemize} 
\item[\underline{2-i.}] 
If $b=a^{-1}$, 
then we have $\phi(a,a^{-1})=\phi(a*c,a^{-1}*c)$. 
By changing $a$ and $c$ variously, 
(2) implies that \underline{$\phi(a,a^{-1})$ is constant} 
regardless of $a\in Q_6$, 
which we denote by $\delta\in A$. 
\item[\underline{2-ii.}] 
If $c=b^{-1}$, 
then the condition (2) is equivalent to 
\begin{eqnarray} 
\phi(a,b^{-1})+\phi(a*b^{-1},b)
-\phi(a,b)-\phi(a*b,b^{-1})=0, 
\end{eqnarray} 
for any pair $(a,b)$ which spans an edge of $H$. 
We also illustrate the condition (4) in Figure~\ref{fig08}. 
\item[\underline{2-iii.}] 
If $c=a^{-1}$, 
then (2) is equivalent to 
\begin{eqnarray}\label{eq05} 
\phi(a,b)+\phi(a*b,a^{-1})-\phi(a,b*a^{-1})=\delta, 
\end{eqnarray} 
for any pair $(a,b)$ which spans an edge of $H$. 
See Figure~\ref{fig08} again, 
where $a*b=b*a^{-1}$ holds. 
\end{itemize} 

\underline{Case 3.} 
Assume that $\{a,b,c\}$ spans 
a face of the octahedron $H$. 
\begin{itemize} 
\item[\underline{3-i.}] 
If $c=b*a$, 
then the condition (2) is equivalent to (5). 
\item[\underline{3-ii.}] 
If $c=a*b$, 
then the condition (2) is equivalent to (3). 
\end{itemize}

We rewrite the elements of $Q_6$ by 
$$1\leftrightarrow(1234), 
2\leftrightarrow(1423), 
3\leftrightarrow(1342), 
4\leftrightarrow(1423), 
5\leftrightarrow(1324), 
6\leftrightarrow(1243).$$ 
Note that 
$1^{-1}=4$, $2^{-1}=5$, and $3^{-1}=6$. 
We consider the map 
$\phi:Q_6\times Q_6 
\rightarrow{\Z}_4=\{0,1,2,3\}$ 
such that $\phi(a,a)=0$ and 
$\phi(a,a^{-1})=1$ for any $a\in Q_6$, and 
\begin{itemize} 
\item $\phi(1,3)=\phi(2,1)=\phi(2,3)=\phi(3,1)=\phi(3,5)$ 

\noindent 
$=\phi(5,1)=\phi(5,6)=\phi(6,1)=\phi(6,2)=\phi(6,5)=1$, 
\item $\phi(1,5)=\phi(5,3)=2$, 
\item $\phi(1,6)=\phi(3,2)=3$, and 
\item $\phi(a,b)=0$ for other cases. 
\end{itemize} 
The value $\phi(a,b)$ 
for $b\ne a,a^{-1}$ is also indicated 
in the lower right of Figure~\ref{fig08} 
by the number of arrows on the edge 
$\overrightarrow{ab}$. 
Then the reader can check that 
$\phi$ satisfies the conditions (3)--(5), 
and hence, $\phi$ is a $2$-cocycle in $Z^2(Q_6;{\Z}_4)$.

\begin{figure}[htb]
\begin{center}
\includegraphics{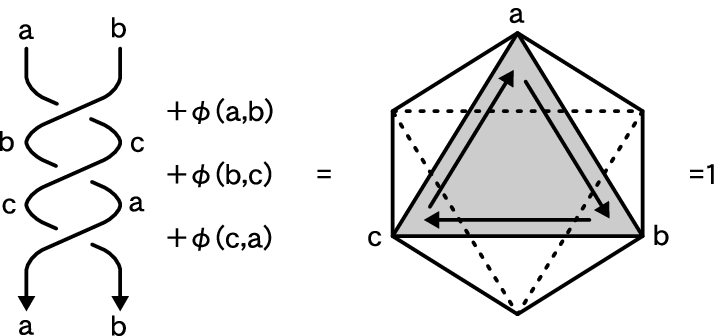}
\caption{}
\label{fig09}
\end{center}
\end{figure}

For this $2$-cocycle $\phi$, 
we calculate the invariant of $T(2,l)$ 
for $l\equiv 3$ (mod $6$). 
Consider the diagram of $T(2,l)$ 
as a closure of the $2$-string braid with $l$ half twists. 
Since each $Q_6$-coloring of the diagram 
is determined by the pair of colors $(a,b)$ 
on the top arcs of the braid, 
we denote the coloring by $C(a,b)$. 
There are $6$ trivial $Q_6$-colorings 
$C(a,a)$ for which we have 
$W_{\phi}\bigl(C(a,a)\bigr)=0$ 
by definition. 
If $b=a^{-1}$, 
then the bottom arcs of the braid 
admits the pair of colors $(a^{-1},a)$; 
for $l$ is odd. 
Hence, such a $Q_6$-coloring does not exist. 
If $b\ne a, a^{-1}$, 
that is, $\{a,b\}$ 
is the boundary of an edge of the octahedron $H$, 
then the same pair of colors appears 
by three half twists. 
See Figure~\ref{fig09}. 
The number of such $Q_6$-colorings are 
$6\times (6-2)=24$. 
For each $Q_6$-coloring $C(a,b)$ with $b\ne a,a^{-1}$, 
we have 
$$W_{\phi}\bigl(C(a,b)\bigr)=
\frac{l}{3}\bigl(\phi(a,b)+\phi(b,c)+\phi(c,a)\bigr),$$
where $c=a*b$. 
On the other hand, 
we see that 
$\phi(a,b)+\phi(b,c)+\phi(c,a)=1$ 
by the definition of $\phi$. 
Hence, we have 
$W_{\phi}\bigl(C(a,b)\bigr)=
\frac{l}{3}\equiv l+2$ (mod $4$), 
and 
\begin{eqnarray*} 
\Phi_{\phi}\bigl(T(2,l)\bigr) &=& 
\bigl\{W_{\phi}\bigl(C(a,b)\bigr)\bigl|
a=b \mbox{ or } b\ne a,a^{-1}\bigr\} \\ 
&=& \bigl\{\underbrace{0,0,\dots,0}_{6}, 
\underbrace{l+2,l+2,\dots,l+2}_{24}\bigr\}.
\end{eqnarray*}

The calculation of $\Phi_{\phi}\bigl(S(m,n)\bigr)$ 
for $m,n\equiv 3$ (mod $6$) 
can be similarly checked, 
and the details are left to the reader. 
The classical knot $S(m,n)$ has a diagram 
as a closure of the $3$-string braid 
$\sigma_1^m\sigma_2^{-1}\sigma_1^n\sigma_2^{-1}$ 
for the standard generators $\sigma_1$ and $\sigma_2$ 
of the braid group. 
Let $C(a,b,c)$ be the $Q_6$-coloring of the diagram 
such that the colors of the top first, second, and third arcs are 
$a,b,c\in Q_6$, respectively. 
Then we have the following three cases; 
\begin{itemize} 
\item $W_{\phi}\bigl(C(a,a,a)\bigr)=0$ 
for any $a\in Q_6$, 
\item $W_{\phi}\bigl(C(a,b,b)\bigr)=m+n+2$ 
for any $a,b\in Q_6$ with $b\ne a,a^{-1}$, and 
\item $W_{\phi}\bigl(C(a,b,b^{-1})\bigr)=m+n$ 
for any $a,b\in Q_6$ with $b\ne a,a^{-1}$. 
\end{itemize} 
Since the numbers of $Q_6$-colorings 
in these cases are $6$, $24$, and $24$, 
respectively, and since $m+n$ is even, 
we have 
$$\Phi_{\phi}\bigl(S(m,n)\bigr)=
\bigl\{\underbrace{0,0,\dots,0}_{6+24=30}, 
\underbrace{2,2,\dots,2}_{24}\bigr\}.$$

\end{document}